\documentclass{article}
\usepackage{amsfonts,amsthm,amsmath,amssymb}

\newtheorem{thm}{Theorem}
\newtheorem{prop}{Proposition}
\newtheorem{lemma}{Lemma}
\newtheorem{cor}{Corollary}
\newtheorem{example}{{\it Example}}
\newtheorem{remark}{Remark}

\def\G{{\rm G}}
\def\dim{{\rm dim}}
\def\Ad{{\rm Ad}} 
\def\GL{{\rm GL}}
\def\deg{\mathop{\rm gdeg}} 
\def\C{\mathbb{C}}
\def\P{\mathbb{P}}
\def\g{\mathfrak{g}}
\def\Oc{{\cal O}}
\def\phi{\varphi}

\def\a{\alpha}
\def\b{\beta}
\def\F{{\cal F}}
\def\S{{\cal S}}
\def\D{{\cal D}}
\def\e{{\rm e}}
\def\Tr{{\rm tr}}
\def\Lie{{\rm Lie}}
\def\T{{\rm T}}
\def\ve{{\varepsilon}}
\def\gl{\mathfrak{gl}}
\def\gdeg{{\rm gdeg}}

\def\Lie{{\rm Lie}}

\title{A Gauss-Bonnet theorem for constructible sheaves on reductive groups}
\author{V. Kiritchenko}
\date{}

\begin{document}
\maketitle
\section {Introduction }  
 In this paper, we prove an analog of the Gauss-Bonnet formula for constructible sheaves on  reductive groups. As a corollary from this formula we get that if a perverse sheaf on a reductive group is equivariant under the adjoint action, then its Euler characteristic is nonnegative.  
    
 In the sequel by a constructible complex we will always mean a bounded complex of sheaves of $\C$-vector spaces whose cohomology sheaves are constructible with respect to some finite algebraic stratification. 

We now formulate the main results. Let $G$ be a complex reductive group, and let $\F$ be a constructible complex on $G$. 
Denote by $CC(\F)$ the characteristic cycle of $\F$. 
It is a linear combination
of Lagrangian subvarieties $CC(\F)=\sum c_\a\T^*_{X_\a} G$ (see \cite{Kash}). Here and in the sequel $\T^*_XG$ denotes the closure of the conormal bundle to the smooth locus of a subvariety $X\subset G$. With $X$ one can associate a nonnegative number $\gdeg(X)$ called the {\em Gaussian degree} of $X$. It is equal to the number of zeros of a generic left-invariant differential 1-form on $G$ restricted on $X$. The precise definition of the Gaussian degree and the Gauss map is given in subsection 2.1. 
 
\begin {thm}\label{main}
If $\F$ is equivariant under the adjoint action of $G$, then its Euler characteristic can be computed in terms of the characteristic cycle by the following formula
$$\chi(G,\F)=\sum c_a \gdeg(X_\a). $$
\end{thm}

 For a perverse sheaf the multiplicities $c_\a$ of its characteristic cycle are nonnegative \cite{Ginz}. The Gaussian degrees of $X_\a$ are also nonnegative by their definition, see below. Thus Theorem \ref{main} immediately implies the following important corollary.

\begin{cor}
If $\F$ is a perverse sheaf equivariant under the adjoint action of $G$,
then its Euler characteristic is nonnegative. 
\end{cor}
In particular, let $\underline\C_X$ be a constant sheaf on a subvariety $X\subset G$ extended by $0$ to $G$. Applying the above statements to this sheaf we get the following corollary.
\begin{cor}\label{cor}
If $X\subset G$ is a closed smooth subvariety invariant under the adjoint action of $G$, then $\chi(X)=(-1)^{\dim X}\gdeg(X)$. Thus the number $(-1)^{\dim X}\chi(X)$ is nonnegative.   
\end{cor}
Indeed, the characteristic cycle of $\underline\C_X $ coincides with $(-1)^{\dim X}\T_X^*G$. Corollary \ref{cor} is a noncompact analog of the classical Hopf theorem which states that the Euler characteristic of a compact oriented $C^\infty$-manifold $M$ is equal to $(-1)^{\dim M}$ times the number of zeros of a generic $1$-form on $M$, counted with signs (coming from the orientation). 

For the case when $G=(\C^*)^n$ is a torus, the Corollary 1 was first proved by F. Loeser and C. Sabbah \cite{LS} with another proof given by  O. Gabber and F. Loeser \cite{GL}. Theorem \ref{main} was  proved in the torus case by J. Franecki
and M. Kapranov \cite{Kapranov}. Theorem \ref{main} holds for all constructible sheaves on a torus. However, it does not hold for arbitrary constructible sheaves on a noncommutative algebraic group (see \cite{Kapranov} for a counterexample). 
M. Kapranov conjectured that it may be still true for constructible sheaves on reductive groups if we consider only sheaves equivariant under the adjoint action. The present paper proves this conjecture. Recently, A. Braverman proved Corollary 1 in a particular case, namely, when a perverse sheaf coincides with its Goreski--MacPherson extension from the set of all regular semisimple elements of $G$ \cite{Brav}.

The main step in the proof of Theorem \ref{main} is to reduce the problem to the case of a maximal torus $T\subset G$. Since $\F$ is $\Ad~G$-equivariant, it is constructible with respect to some Whitney stratification $\S$ with $\Ad~ G$-invariant strata. In section \ref{s.eu} we prove that the Euler characteristic of a stratum $X\in \S$ coincides with that of the intersection $X\cap T$. This implies that the sheaf $\F$ restricted onto the maximal torus $T$ has the same Euler characteristic as $\F$. In subsections $2.2,2.3$ we recall some  facts about Euler characteristic needed for the proof.
In section \ref{gauss} we prove that the Gaussian degrees of $X$ and $X\cap T$ coincide. 

To deal with the characteristic cycle we use the Dubson-Kashiwara index formula that expresses the multiplicities $c_a$ in terms of the local Euler characteristic of $\F$ along 
each stratum and some topological data depending on the stratification only (subsection 2.4). This data is given by the Euler characteristics with compact support of complex links. In our case we can choose a complex link to be invariant under the action of some torus and thus simplify computation of its Euler characteristic. This approach is taken from \cite{Mirc}.
In section \ref{s.li} we prove that for any stratum $X_\b\in\S$ and any semisimple stratum $X_\a\in\S$, such that $X_\a\subset\overline{X_\b}$, the  Euler characteristic with compact support of their complex link  coincides with that of the complex link of the strata $X_\a\cap T$ and $ X_\b\cap T$ . This allows us to view the formula from Theorem \ref{main} as the same formula for the restriction of $\F$ onto $T$
(section \ref{proof}). Then we apply the result of \cite{Kapranov}.

I am grateful to M. Kapranov for constant encouragement, valuable suggestions and remarks. 
I would also like to thank A. Braverman  and A. Khovanskii for useful discussions.

\section{Preliminaries} 

\subsection{Gaussian degree}
 We now define the (left) Gauss map and the Gaussian degree. The material of this subsection is taken from \cite{Kapranov}. For more details see \cite{Kapranov,Fra}.

 Let $G$ be a complex algebraic group with Lie algebra $\g$, and let $X$ be its subvariety
of the dimension $k$. Denote by $\G(k,\g)$ the Grassmanian of $k$-dimensional subspaces in $\g$. For any point $x\in G$ there is a natural isomorphism between the tangent space $\T_xG$ and $\g$ given by the left multiplication by $x^{-1}$: 
$$ L_x: y\mapsto x^{-1}y; \quad d_xL_x: \T_xG\to\g.$$
The left Gauss map $\Gamma_X: X\to \G(k,\g )$ is defined as follows:
$$\Gamma_X(x)=d_xL_x(\T_xX).$$
The Gauss map is rational and regular on the smooth locus $X^{sm}$ of $X$.  
If $X$ is a hypersurface, $\Gamma_X$ maps $X$ to $\P(\g^*)$, which has the same dimension as $X$. In this case we define the Gaussian degree of $X$ as the degree of its Gauss map. In general case the Gaussian degree is the degree of a rational map $\tilde \Gamma_X:\tilde X\to\P(\g^*)$, where $\tilde X$ and $\tilde\Gamma_X$ are defined as follows. The variety $\tilde X$ is a fiber bundle over $X^{sm}$, whose fiber at a point $X$ consists of all hyperplanes in $\g$ that contain a subspace $\Gamma_X(x)$, i.e.
$\tilde X=\{(x,y)\in X^{sm}\times \P(\g^*): \Gamma_X(x)\subset y\}$. Then $\tilde\Gamma_X(x,y)=y$. Note that $\tilde X$ and $\P(\g^*)$ have the same dimension. It is clear from the definition that the Gaussian degree is a birational invariant of a subvariety.

In the sequel we will use another description of the Gaussian degree. Let $\omega$ be a  generic  left-invariant differential 1-form on $G$ ( for reductive groups we define a generic 1-form explicitly in section \ref{gauss}). We call a point $x\in X$ a {\em zero} of $\omega$, if $\omega$ restricted to the tangent space $\T_xX$ is zero. Then it is easy to verify that the Gaussian degree of $X$ is equal to the number of zeros of $\omega$ on $X$. 

\subsection{Euler characteristic} 

Let $T$ be a torus ($T$ may be a complex torus $(\C^*)^n$ as well as a compact one $(S^1)^n$). Consider its linear algebraic action on $\C^N$, and a locally closed semialgebraic subset $X\subset\C^N$  
invariant under this action. Let $X^T\subset X$ be the set of the fixed points. In what follows $\chi$ denotes the usual topological Euler characteristic and $\chi^c$ the Euler characteristic computed using cohomology with compact support. The following simple and well-known fact plays the crucial role in the sequel. 

\begin{prop}\label{tor} 
The spaces $X$ and $X^T$ have the same Euler characteristic with compact support: 
$$\chi^c(X)=\chi^c(X^T).\eqno\square$$
\end{prop}

\subsection{The Euler characteristic of sheaves}
  We now recall a formula for the Euler characteristic of constructible sheaves on varieties.
 Let $X\subset\C^N$ be a smooth subvariety, and let $\F$ be a constructible complex on $X$.  With any point $x\in X$ we can associate the local Euler characteristic $\chi(\F_x)$ of $\F$ at this point. Thus $\F$ gives rise to the constructible function $\chi(\F)$ on $X$ by the formula $\chi(\F)(x)=\chi( F_x)$. 

There is the 
concept of the direct image of a constructible function (see \cite{Fulton} and \cite{Kash}, Section 9.7). It is defined for any morphism of algebraic varieties 
$X\to Y$ and a constructible function on $X$. We use the more suggestive notation $\int_X f(x)d\chi$ for the direct image of $f$ under the morphism $X\to pt$, since this direct image may be also defined as the integral of $f$ over the Euler characteristic \cite{Viro,Askold}. 

\begin{prop}\label{const}
The global Euler characteristic $\chi(X,\F)=\sum (-1)^iH^i(X,\F)$  is equal to the following integral over the Euler characteristic
$$\chi(X,\F)=\int_X \chi(\F)d\chi$$
In other words, if we fix an algebraic stratification $X=\bigsqcup X_\a, \a\in\S$ , such that the function $\chi(\F)$ is constant along each stratum, we get
$$\chi(X,\F)=\sum_{\a\in\S}\chi_\a(\F)\chi^c(X_\a),$$
where $\chi_{\a}(\F)$ is the value of $\chi(\F)$ at any point of a stratum $X_\a$.
\end{prop}  
See \cite{Kash}, Section 9.7 for the proof.

\subsection{Complex links and characteristic cycles}
 Let us recall the definition of a complex link.
 We use the notations of the previous subsection.
 Suppose that $\S$ is a Whitney stratification of $X$.
Let $X_\a,X_\b,\a,\b\in\S$, be two strata such that $X_\a\subset\overline{X_\b}$. Choose a point $a\in X_\a$ and any normal slice $N\subset X$ to $X_\a$ at the point $a$.
Let $v^*$ be a generic covector in the fiber of the conormal bundle $\T^*_{X_\a}\C^N$ at the point $a$ (i.e. $v^*$ belongs to some open dense subset of this fiber that depends on the stratification $\S$). This covector defines a linear function $l$ on $N$ by the formula $l(x)=v^*(x-a)$. Let $h(.,.)$ be a Hermitian metric in $\C^N$ and $B=\{x\in \C^N: h(x-a,x-a)\le const\}$ a small ball with the center at $a$. We now define the complex link $L$ of the strata $X_\a, X_\b$ as
$L=B\cap l^{-1}(\varepsilon)\cap X_\b$. If $\varepsilon$ and the radius of the ball $B$ are small enough, the result up to a homeomorphism does not depend on any of the choices involved. This is proven in \cite{MacPhG}. We will use the notation $\e(\a,\b)$ as well as $\e(X_\a,X_\b)$ for 
the Euler characteristic of $L$ with compact support. We also set $\e(\a,\a)=-1$.   

The numbers $\e(\a,\b)$ are useful when one need to find the multiplicities of the characteristic cycle $CC(\F)$. Multiplicities are recovered from the constructible function $\chi(\F)$  by the following theorem of Dubson and Kashiwara.

\begin{thm}\label{mult}
The characteristic cycle of $\F$ is the linear combination of Lagrangian subvarieties $\T^*_{X_\a}X, \a\in\S$ with coefficients
$$c_\a=(-1)^{\dim X_\a +1} \sum_{X_\a\subset\overline{X_\b}}\e(\a,\b)\chi_\b(\F).$$ 
\end{thm}
See \cite{Ginz}, Theorem 8.2 for the proof.

\section{Euler characteristic of invariant subvarieties}

\label{s.eu}
Let $G$ be a connected reductive group over $\C$, and $T$ a maximal complex torus in $G$. Consider a subvariety $X\subset G$ invariant under the adjoint action of $G$.

\begin{prop}
\label{euler}
The Euler characteristics with compact support of $X$ and of $X\cap T$ coincide.
If in addition $X$ and $X\cap T$ are smooth, then $\chi(X)=\chi(X\cap T)$.
\end{prop}

\begin{proof}
The subvariety $X$ is invariant under the adjoint action of $G$. In particular, it is invariant under the adjoint
action of the maximal torus $T$. The set $G^T\subset G$ of the fixed points under this action 
coincides with $T$, since the centralizer of the maximal torus coincides with the maximal torus 
itself. Thus by Proposition \ref{tor} the varieties $X$ and $X\cap T$ have the same Euler characteristic with compact support. 
If $X$ and $X\cap T$ are smooth, then by Poincar\'e duality their Euler characteristics with compact support coincide with the usual ones, and $\chi(X)=\chi(X\cap T).$
\end{proof}

\begin{example} \em Let $X=\Oc_a$ be the orbit of an element $a\in G$ under the adjoint action of $G$.
Then proposition \ref{euler} implies that if $a$ is semisimple, then $\chi(\Oc_a)$ is equal to 
the number $|\Oc_a\cap T|$ of intersection points. We may choose the maximal torus $T$, such that
$a\in T$. Since the orbit of $a$ under the action of the Weyl group $W$ on $T$ coincides with 
$\Oc_a\cap T$, we obtain that $\chi(\Oc_a)=|W|/|Stab~a|$, where $Stab~a\subset W$ is the stabilizer
of $a$ in $W$. If $a$ is not semisimple, then $\chi(\Oc_a)=0$.
\end{example}

Let $\F$ be a constructible complex on $G$.

\begin{prop}\label{prop}

Suppose that $\F$ is equivariant under the adjoint action of $G$. Let $\F_T$ be a restriction of $\F$ onto $T\subset G$. Then the sheaves $\F$ and $\F_T$ have the same Euler characteristic:
$$\chi(G,\F)=\chi(T,\F_T).$$
\end{prop}
\begin {proof}  The sheaves $\F$ and $\F_T$ have the same local Euler characteristic at a point $x\in T$, since $\F_T$ is the restriction of $\F$ on $T$. Thus the Euler characteristic $\chi(T,\F_T)$ is equal to $\int_T \chi(\F)d\chi$ by Proposition \ref{const}. The function $\chi(\F)$ is invariant under the adjoint action of $G$, and Proposition \ref{euler} implies $$\int_T \chi(\F)d\chi=\int_G \chi(\F)d\chi.$$
The last integral is equal to the Euler characteristic $\chi(G,\F)$ by 
Proposition~\ref{const}.
\end {proof}

\section{Gaussian degree of invariant subvarieties}\label{gauss}
We now compare the Gaussian degrees of $X$ in $G$ and of $X\cap T$ in $T$. 
 Clearly, the Gaussian degree of a $k$-dimensional subvariety is equal to the sum of the Gaussian degrees of its $k$-dimensional irreducible components. Thus we may assume that $X$ is irreducible. There are two 
cases: the set of all nonsemisimple elements of $X$ has codimension either 0 or at least 1. In what follows we prove that in the first case $\gdeg(X)=0$, and in the second case 
$\gdeg(X)=\gdeg(X\cap T)$. In particular, the Gaussian degree of any orbit $\Oc_a\subset G$ 
coincides with the Gaussian degree of its intersection with the maximal torus $T$. 

Any reductive group $G$ admits an embedding in $\GL_N(\C)$ for some $N$, such that the inner product $\Tr(Y_1Y_2)$ is nondegenerate on Lie algebra $\g$. Let us fix such an embedding. 
Then $\g$ may be identified with the space of all left-invariant differential $1$-forms on $G$: an element $S\in\g$ gives rise to a 1-form $\omega$ by the formula 
$$\omega(Y)=\Tr(x^{-1}YS), \eqno(1)$$ 
where $x\in G$ and $Y\in\T_xG$. We will call such a form {\em generic}, if $S$ is regular semisimple.

\begin{lemma}
 All generic left invariant $1$-forms form an open dense subset in the space of all left-invariant $1$-forms.
\end{lemma}
\begin{proof}
All regular semisimple elements  form an open dense subset in $\g$.
This implies the statement of the lemma.
\end{proof}
\begin{prop}
\label{orbit}
The Gaussian degree of an orbit $\Oc_a$ is equal to the number of intersection points
$\Oc_a\cap T$. In particular, if $a$ is a nonsemisimple element, then $\deg(\Oc_a)=0.$
\end{prop}

\begin{proof}
Consider the map
$$\phi: G\to\Oc_a;\quad \phi: g\mapsto gag^{-1}.$$
Since $\phi$ is smooth and surjective, the tangent space $\T_x\Oc_a$ is the image of the induced 
map $d\phi$. A simple computation shows that $\T_x\Oc_a=[\g,x]$. Let $\omega$ be a generic left-invariant differential 1-form
on $G$ given by the formula (1). Then $\omega=0$ on $\T_x\Oc_a$ is equivalent to 
$\Tr(x^{-1}YxS-YS)=0$ for any $Y\in\g$. Since the form $\Tr(Y_1Y_2)$ is $\Ad~G$-invariant,
we have $\Tr(x^{-1}YxS-YS)=\Tr(Y(xSx^{-1}-S))$. The form $\Tr(Y_1Y_2)$ is nondegenerate on $\g$. Thus $x$ and $S$ commute, and $x$ belongs to some maximal torus $T_S$ that depends on $S$. 
\end{proof}

\begin{remark}
Suppose that an element $a$ lies in the maximal torus $T$. Then the space $\T_a\Oc_a=[\g,a]$ is orthogonal to the tangent space $\T_aT$ with respect to the form $(Y_1,Y_2)\mapsto\Tr(a^{-1}Y_1\cdot a^{-1}Y_2)$. Since this form is nondegenerate on $\T_aT$, we get that at any point $x\in\Oc_a\cap T$ the intersection of tangent spaces $\T_x\Oc_a$ and $\T_xT$  is zero.
\end{remark}  

\begin{cor}
\label{nilp}
Let $Z\subset G$ be an irreducible subvariety invariant under the adjoint action of $G$, such that the 
set $Z_n$ of all nonsemisimple elements of 
$Z$ is a Zariski open nonempty subset in $Z$. Then $\deg(Z)=0$. 
\end{cor}
\begin{proof}
The Gaussian degree of a subvariety $Z\subset G$ is birationally invariant. Therefore, it  
suffices to compute it for $Z_n$. Let $\omega$ be a generic left-invariant differential 1-form
on $G$ given by the formula (1).  For any smooth point $a\in Z_n$  the restriction of this form on the 
subspace $\T_a\Oc_a\subset\T_aZ_n$ of the tangent space $\T_aZ_n$ is already nonzero by 
Proposition \ref{orbit}. Thus the form $\omega$ does not vanish in any smooth point of $Z_n$, 
and $\gdeg(Z)=\gdeg(Z_n)=0$. 
\end{proof}

\begin{prop}
\label{semis}
Let X be an irreducible subvariety of $G$ invariant under the adjoint action of $G$, such that the 
subset of all nonsemisimple elements of $X$
has codimension at least $1$ in $X$. Then the Gaussian degrees of $X$ in $G$ and $X\cap T$ in $T$ 
coincide. 
\end{prop}

\begin{proof}
Let $k$ be the maximal dimension of a semisimple
orbit in X. Denote by $X_s$ the set of all semisimple 
elements in $X$, whose orbits have dimension $k$. Then $X_s$ is a Zariski open subset of $X$.
Consider the map 
$$\phi:G\times(X\cap T)\to X,\quad (g,t)\mapsto gtg^{-1}.$$
The image of $\phi$ contains $X_s$. Since $X_s$ is a Zariski open nonempty subset of $X$, 
$\deg(X)=\deg(X_s)$. For any smooth point $x\in X_s$ the tangent space $\T_xX$ is again the image of the 
induced map $$d\phi:\T_gG\times \T_t(X\cap T)\to \T_xX, $$
where $gtg^{-1}=x$. By calculating $d\phi$ we obtain that 
$\T_xX=[\g,x]\oplus g\T_t(X\cap T)g^{-1}$. 
Let $\omega$ be a generic left-invariant differential 1-form
on $G$ given by the formula (1). 
Then $\omega=0$ on $\T_xX$ is equivalent to
$\omega=0$ on $[\g,x]$ and $\omega=0$ on $g\T_t(X\cap T)g^{-1}$. The first identity holds if and 
only if $x$ belongs to the maximal torus $T_S$ (see the proof of Proposition \ref{orbit}). Denote by $\omega_T$ 
the restriction of $\omega$ on $T^*T_S$. If $x\in T_S$, then 
$g\T_t(X\cap T)g^{-1}=\T_x(X\cap T_S)$. Thus the form $\omega$ vanishes on $\T_xX$ if and only if 
the form $\omega_T$ vanishes on $\T_x(X\cap T_S)$. It follows that $\deg(X)=\deg(X\cap T_S)=\deg
(X\cap T)$, since all maximal tori are conjugate. 
\end{proof}

\section{The Euler characteristic of the complex link}
\label{s.li}
We now compute the Euler characteristic with compact support of a complex link for a certain class of stratifications of $G$. 
For any $a\in G$ we define the rank of $a$ to be the dimension of its centralizer in $G$. 
A Whitney stratification $\S$ of $G$ is called {\em admissible} if the following conditions hold. For every $\a\in\S$

\begin{itemize}
\item the stratum $X_\a$ is invariant under the adjoint action of $G$,
\item elements of $X_\a$ are either all semisimple of the fixed rank or all nonsemisimple. 
\end{itemize}

Due to the second condition and Remark 1,  for any semisimple stratum $X_\a\in S$ intersection $X_\a\cap T$ is smooth, and at any point $x\in X_\a\cap T$ the intersection $\T_x{X_\a}\cap \T_xT$ coincides with $\T_x(X_\a\cap T)$. Thus we can consider an induced Whitney stratification $\S_T$ of the maximal torus $T$, namely, 
$T=\bigsqcup (X_\a\cap T), \a\in\S_0,$ where $\S_0\subset\S$ consists of all semisimple strata.

\begin {prop}
\label{link}
Consider two strata $X_\a$ and $X_\b$, such that $X_\a$ belongs to the closure of $X_\b$.
If $X_\a$ is semisimple and $X_\b$ is not, then $\e(\a,\b)=0$. 
If both $X_\a$, $X_\b$ are semisimple, then $\e(\a,\b)=\e(X_\a\cap T,X_\b\cap T)$, where 
the complex link of $X_\a\cap T$ and $X_\b\cap T$ is taken in the torus $T$. 
\end {prop}

\begin{proof}
Let $Z\subset G$ be the centralizer of an element $a\in X_\a$. Then $Z$ is again a reductive group. Since the tangent spaces $\T_a Z$ and $\T_a \Oc_a$ are orthogonal with respect to the form $(Y_1,Y_2)\mapsto\Tr(a^{-1}Y_1\cdot a^{-1}Y_2)$, and this form is nondegenerate on $\T_a Z$, we get that $Z$ is the normal slice to the orbit 
$\Oc_a\subset X_\a$. Thus any normal slice to $Z\cap X_\a$ in $Z$ will also be the normal slice to $X_\a$ in $G$. Let us construct a normal slice $N\subset Z$ invariant under 
the adjoint action of $Z$. 
 
Let $k$ be a dimension of $X_\a\cap Z$.
 Some neighborhood of $a$ in $X_\a\cap Z$ lies in
the center of $Z$, because all elements of $X_\a$ have
the same rank.  Thus we can find $k$ characters 
$\phi_1,\ldots,\phi_k$ of the group $Z$, such that their differentials $d_a\phi_1,\ldots,d_a\phi_k$
restricted to the tangent space $\T_a(X_\a\cap Z)$ are linearly independent. 
Let $N\subset Z$ be the set of common zeros of the system $\phi_1(za^{-1})=\ldots=\phi_k(za^{-1})=1$. 
 
\begin{example}\em 
a) Let $G$ be $\GL_N(\C)$ and let $X_\a=Z(\GL_N)=\C^*$ be the center of $\GL_N$. Then $Z=G$ and the only characters of $Z$ are the powers of determinant. We have $d_edet=\Tr$ for the identity element $e\in\GL_N$, and $\Tr$ is nonzero linear function on $\C^*e$. Thus at the point $e\in X_\a$ we can take $N={\rm SL}_N$.

b) Let $G$ be any reductive group, and let $X_\a$ be a stratum consisting of regular semisimple elements. Then $Z$ is a maximal torus. Thus any normal slice to $X_\a\cap T$ in $Z$ is invariant under the adjoint action of $Z$.
\end{example}

 We now continue the proof of Proposition \ref{link}. Consider a generic linear function $l(x)$ 
on $N$ given by the formula $l(x)=\Tr((a^{-1}x-e)S)$, where $S\in\Lie~Z$ is 
regular semisimple. 
There exists a maximal torus $T\subset Z$ centralizing
$S$. Since $a$ is semisimple, $T$ is also the maximal torus in $G$. For any $\ve$ the hyperplane 
$l^{-1}(\ve)$ is invariant under the adjoint action of $T$. 
Denote by $T_c$ the compact form of $T$. Choose a Hermitian inner product $h(.,.)$ on 
$\gl_N$ invariant 
under the adjoint action of $T_c$ and a small ball $B=\{x\in\gl_n: h(x-a,x-a)\le const\}$. 

Thus with a generic vector $S\in \Lie Z$ we associate the complex link 
$L=B\cap l^{-1}(\ve)\cap X_\b$ 
of the strata $X_\a$ and $X_\b$. The complex link $L$ is invariant under the adjoint action of 
the torus $T_c$ by the construction. 
Thus by Proposition \ref{euler} we get $\chi^c(L)=\chi^c(L^{T_c})=\chi^c(L\cap Z^{T_c})$. Note that $Z^{T_c}=Z^T=T$. 
If $X_\b$ is nonsemisimple, $X_\b\cap T$ is empty, thus $L\cap T$ is empty. It follows that $\e(\a,\b)=0$ in this case. If $X_\b$ is semisimple, 
then $L\cap T$ is a complex link for $X_\a\cap T, X_\b\cap T$ in the torus $T$.        
\end{proof}
 
\begin{cor}
If $X$ is a smooth irreducible subvariety invariant under the adjoint action of $G$, then either $X$ consists of nonsemisimple elements only or the set of all semisimple elements in  $X$ is dense. In particular, if in addition $X$ is closed, then it contains a dense subset of semisimple elements. The Gaussian degrees of $X$ and $X\cap T$ coincide.
\end{cor}  
\begin{proof}
Let us prove the first statement by contradiction. Let $\S$ be an admissible stratification of $G$ subordinate to $X$, and let $X_n\subset X$ be a maximal open stratum in $X$, such that $X_n$ is nonsemisimple. Then $\dim X_n=\dim X$, and $X-X_n$ contains at least one semisimple stratum $X_s$. The number $\e(X_s,X_n)$ is zero by Proposition \ref{link}. That contradicts to the smoothness of $X$. Now combining the first statement with the results of Section \ref{gauss} we get the last statement.
\end{proof}

\section{Proof of Theorem \ref{main}} 
 \label{proof}
 Since $\F$ is constructible and equivariant under the adjoint action, there exists some finite algebraic Whitney stratification $\S$ subordinate to $\F$, such that each stratum is invariant under the adjoint action. Stratifying each stratum if necessary we may assume that $\S$ is admissible. 
Denote by $i:T\to G$ the embedding of a torus $T\subset G$. 
Let us apply Theorem \ref{mult} and Proposition \ref{link} to the characteristic cycles of $\F$ and $\F_T$. Notice that for a semisimple stratum $X_\a\in\S$ the difference $\dim X_\a-\dim (X_\a\cap\S)$ is equal to $\dim \Oc_a, a\in X_\a$, and the latter is even. As a straightforward corollary we get that the inverse image $i^*CC(\F)$  coincide with $CC(\F_T)$ on the set of all semisimple strata. More precisely, the following is true.  
  
\begin{cor}\label{cycle}
Let $X_\a\in\S$ be a semisimple stratum. The multiplicities of characteristic cycles of $\F$
and $\F_T$ along the strata $X_\a$ and $X_\a\cap T$, respectively, coincide. 
\end{cor}
 
\begin{example} \em
Suppose that the support of the constructible function $\chi(\F)$ lies in the closure of an orbit $\Oc_a, a\in G$. This kind of sheaves is studied in \cite{Mirc} for unipotent orbits. In this case the strata of an admissible stratification that contribute to the the characteristic cycle are the orbits in $\overline{\Oc_a}$. Let $a_s,a_n\in G$
 be the semisimple and nilpotent elements respectively, such that $a=a_s+a_n$. Then
$X_\a=\Oc_{\a_s}$ is the only semisimple stratum in $\overline{\Oc_a}$. Thus for the multiplicity $c_\a(\F)$ of $CC(\F)$ along this stratum we get $c_{\a}(\F)=\chi_\a(\F)=c_{\a}(\F_T)$.  
\end{example}

Now the formula of Theorem \ref{main} reduces to the same formula for the sheaf $\F_T$ and the stratification $\S_T$.
First, $\chi(\F,G)=\chi(\F_T,T)$ by Proposition \ref{prop}. Second, for all nonsemisimple strata $X_\a,\a\in\S$ we have $\gdeg(X_\a)=0$ by Proposition \ref{nilp}. Thus the right hand side of the formula may be considered as the sum over semisimple strata only, i.e.
$$\chi(X,\F)=\sum_{\a\in\S_0}c_\a(\F)\gdeg(X_\a).$$
By Corollary \ref{cycle} this is equivalent to the formula 
$$\chi(T,\F_T)=\sum_{\a\in\S_T}c_\a(\F_T)\gdeg(X_\a\cap T),$$ 
since $\gdeg(X_a)=\gdeg(X_a\cap T)$ by Proposition \ref{semis}. To prove the latter formula we apply the Theorem 1.3 from \cite{Kapranov}.

\end{document}